\RequirePackage{fix-cm}
\documentclass[smallextended]{svjour3}       
\smartqed  
\usepackage{graphicx}
\usepackage{amssymb}
\usepackage{amsmath} 
\usepackage[linesnumbered,ruled]{algorithm2e} 
\usepackage{hyperref}  

\usepackage{algorithmic}
\usepackage{xcolor}
\usepackage{subfigure}  

\DeclareMathOperator*{\scup}{\sqcup}

\makeatletter
\newcommand{\labitem}[2]{%
	\def\@itemlabel{\textsc{#1}}
	\item
	\def\@currentlabel{\textsc{#1}}\label{#2}}
\makeatother

\makeatletter
\newcommand{\labitemnr}[2]{%
	\def\@itemlabel{#1}
	\item
	\def\@currentlabel{#1}\label{#2}}
\makeatother

\newcommand{\mltg}{\hat{G}}
\newcommand{\mltv}{\hat{N}}
\newcommand{\mlte}{\hat{E}}
\newcommand{\mltc}{\hat{C}}
\newcommand{\usg}{\tilde{G}}
\newcommand{\usv}{\tilde{N}}
\newcommand{\use}{\tilde{E}}

\newcommand{\usy}{\tilde{Y}}
\newcommand{\mltn}{\hat{\nu}}
\begin{document}

\title{Online Barycenter Estimation of Large Weighted Graphs. A multiscale strategy. 
}

\titlerunning{Online Barycenter Estimation of Large Weighted Graphs.}        

\author{Ioana Gavra         \and
        Laurent Risser 
}

\authorrunning{Gavra I. and Risser L.} 

\institute{Ioana Gavra \at
              Toulouse School of Economics.\\
              \email{Ioana.Gavra@math.univ-toulouse.fr}           
           \and
           Laurent Risser \at
              CNRS, Institut de Math\'ematiques de Toulouse.\\
              \email{lrisser@math.univ-toulouse.fr}           
}


\maketitle

\begin{abstract}
In this paper, we propose a new method to compute the barycenter of large weighted graphs endowed with probability measures on their nodes. We suppose that the edge weights are distances between the nodes and that the probability measure on the nodes is related to events observed there. For instance, a graph can represent a subway network: its edge weights are the distance between two stations, and the observed events at each node are the subway users getting in or leaving the subway network at this station. The probability measure on the nodes does not need to be explicitly known. Our strategy only uses observed  node related events to give more or less emphasis to the different nodes. Furthermore, the barycenter estimation can be updated in real time with each new event. 
We propose a multiscale extension of \cite{us} where the decribed strategy is valid only for medium-sized graphs due to memory costs. Our multiscale approach is inspired from the geometrical decomposition of the barycenter in a Euclidean space: we apply a heuristic \textit{divide et impera} strategy based on a preliminary clustering. 
Our strategy is finally assessed on road- and social-networks of up to $10^6$ nodes. We show that its results compare well with \cite{us} in terms of accuracy and stability on small graphs, and that it can additionally be used on large graphs even on standard laptops.

\keywords{Large graphs analysis \and Barycenter estimation \and Online statistics}
\end{abstract}

\section{Introduction}
\label{intro}

\subsection{Context}

\paragraph*{Why graphs? $\mbox{ }$} 
Graph structures can model complex phenomena of high interest in a wide variety of domains and play an important role in various fields of data analysis. Although graphs  have been used for quite a while in some fields, \textit{e.g.}  in sociology since the 1930's \cite{moreno}, the recent explosion of available data and computational resources boosted the importance of studying and understanding networks. Among the main  application fields, one can count computer science (Web understanding \cite{web} and representation), biology (neural or protein networks, genes), social sciences (analysis of citations graphs, social networks \cite{social}), machine learning \cite{Goldenberg}, statistical or quantum physics \cite{estrada}, marketing (consumers preference graphs) and computational linguistics \cite{semg}.

\paragraph*{Barycenter: motivation and applications. $\mbox{ }$}
Singling out the most influential node or nodes can be seen as a first step to understand the structure of a network. Different notions of \emph{node centrality} have been introduced to measure the influence or the importance of nodes of interest in a network. 
Centrality notions are sometimes related to the mean distance from each node to all others \cite{meandist}, to the degree of each node \cite{degree} or even to the eigenvalues of the graph's adjacency matrix  \cite{eigenvalue}. A rather complete survey can be found in \cite{borgatti}. \\
As far as the authours know, these notions of centrality do not take into account any weight on the nodes (but only on the edges), although there are numerous applications where this would be rather natural. For example, in the case of a metro network, when trying to establish a \emph{central} station, it is quite reasonable to take into account the number of passengers that use each station. In the case of a traffic network, the node-weight can model how many cars pass by a given intersection; in the case of a social network it can model the number of followers (or likes, or posts, \textit{etc.}) of each individual. \\ 
To take this kind of information into account, throughout this paper, we interest ourselves to the barycenter of a graph with respect to a probability measure on the node set, as defined in our previous work \cite{us}. As we will see later on, is a natural extension of the expected value on a Euclidean space. Furthermore, our algorithm is developed in an online context: it does not need the exact knowledge of the probability measure (the number of passengers that use each station), but only observations of this random distribution (we can see when a passenger uses a station), and can  be easily updated at the arrival of a new observation.


Besides determining a central node, the knowledge of such a barycenter on a graph can be of multiple use.  For example, from a statistical point of view, for a fixed graph, the computation of the barycenter using two data sets of observations could be used to determine if the two sets are sampled from the same probability measure (on the nodes set). 
Such a mean position can also be a preliminary step for a more detailed study, like the one provided by a generalization of a Principal Component Analysis, that could translate the main statistical fluctuations among the nodes of the network. 
The barycenter can also be useful in graph representation, since setting the barycenter in a central position can provide an intuitive visualization.

\subsection{Online graph barycenter estimation}

\subsubsection{Graph barycenter definition based on  Fr\'echet means}

Since the networks studied in this paper in are finite, 
the node weights can be seen as a probability measure on the nodes set.
For a probability measure defined on an Euclidean space, there are two classical notions
of centrality: the median and the Euclidean mean. Defining an average or central position in a non-euclidean metric space is not straight forward since the natural addition or averaging operations are not necessarily defined.   

Back in $1948$, M. Fr\'echet presented a possible answer to this problem, not only for the median and the mean of a probability measure, but for moments of all orders  \cite{frechet}. He introduced a notion of typical position of order $p$ for a random variable $Z$ defined on a general metric space  $(\mathcal{E},d)$ and distributed according to any probability measure $\nu$. This is now known as the $p$-Fr\'echet mean, or simply the $p$-mean, and is defined as:

$$
M_{\nu}^{(p)} := \arg \min_{x \in \mathcal{E}} \mathbb{E}_{Z \sim \nu}[d^p(x,Z)].
$$

This definition might seem counter-inutitive, but one can notice that this variational formulation also holds for real random variables. For example, if $Z$ is a random variable distributed according to a distribution $\nu$ on $\mathbb{R}^d$, its expected value, given by
 $
 m_{\nu} = \int_{\mathbb{R}^d} x d\nu(x)
 $ is also the point that minimizes:
$$
x \longmapsto \mathbb{E}_{Z \sim \nu} [ |x-Z|^2].
$$

Now, let $G=(N,E)$ denote a finite weighted graph, $E$ its edges set and $\nu$ a probability measure on its nodes set $N$.
The barycenter of a graph $G=(N,E)$ is then naturally defined as the $2$-Fr\'echet mean, that we simply denote Fr\'echet mean:
$$M_{\nu}=\mathrm{argmin}_{x\in N}\sum_{y\in N} d^2(x,y) \nu(y).$$

\subsubsection{Online estimation framework}
We place ourselves in the  online estimation framework, in the sense that we suppose that the probability measure $\nu$ unknown. A sequence $(Y_n)_{n\ge 0}$ of i.i.d. random variables distributed according to $\nu$ is instead available. For instance, an observation $Y_n$ can be interpreted as the access of a passenger to a given station for subway networks, the passage of a car on a given crossroad for traffic networks, or a paper download in a scientific social network.

\section{Barycenter estimation using simulated annealing}\label{sec:SA}

The authors proposed in \cite{us} a method to estimate the barycenter of weighted graphs, based on a simulated annealing algorithm with homogenization. 
In addition to introduce this method, 
we also established in \cite{us} its convergence from a theoretical point of view. Since this method is one of the corner stones of our current work, we briefly explain its principles and its main parameters in this section. We also give  in Alg.~\ref{algo:HLSA} a simplified  pseudo-code that explains how it practically works.

\begin{algorithm}
\caption{Graph barycenter estimation algorithm of \cite{us}}
\label{algo:HLSA}
\begin{algorithmic}[1]
  \REQUIRE Continous version of $G=(N,E)$, \textit{i.e.} $\Gamma_G$.
  \REQUIRE Observations sequence $Y=(Y_k)_{k\ge 1}$ on the nodes set $N$.
  \REQUIRE Increasing inverse temperature $(\beta_t)_{t \geq 0}$ and intensity $(\alpha_t)_{t \geq 0}$.
  \STATE Pick $X_0 \in \Gamma_G$ and set $K=\mathrm{len}(Y)-1$.
  \STATE $T_0=0$.
  \FOR{$k=0:K$}
     \STATE Generate $T_k$ according to $\alpha_k$.
     \STATE Generate $\varepsilon_k \sim \mathcal{N}(0,\sqrt{T_k-T_{k-1}})$.
     \STATE Randomly move $X_k$ (Brownian motion): $X_k=X_k+\bf{h_k}\varepsilon_k$, where $\bf{h_k}$ is a direction uniformly chosen among the directions departing from $X_k$, and $\varepsilon_k$ is a step size.
     \STATE Deterministically move $X_k$ towards $Y_{k+1}$: $X_{k+1}=X_k + \beta_{T_k} \alpha_{T_k}^{-1} \bf{X_kY_{k+1}}$, where $\bf{X_k Y_{k+1}}$ represents the shortest (geodesic) path from $X_k$ to $Y_{k+1}$ in $\Gamma_G$.
  \ENDFOR
  \RETURN Graph location $X_K$ estimated as the barycenter of $\Gamma_G$. We consider the nearest node to $X_K$ in $G$ as its barycenter.
\end{algorithmic}
\end{algorithm}

Simulated annealing is an optimization technique based on a gradient descent dynamic to which we add a random perturbation in order to help it escape local traps. The importance of this random perturbation is then decreased progressively in order to cool down the system  and let the algorithm converge (or stabilize). This effect is parametrized by a continuous function $(\beta_t)_{t\ge 0}$, that represents the inverse of the so-called temperature schedule: when $\beta_t$ is small, the system is \textit{hot} and the random noise is quite important with respect to the gradient descent term. Then, when $\beta_t$ goes to infinity, the random perturbation is negligible. 

Another important parameter comes from the on-line aspect of the algorithm. 
In our model, we simulate the arrival times $T_n$ of the observations $Y_n$ by an inhomogeneous Poisson process $(N^{\alpha}_t)_{t\ge 0}$ \footnote{$T_n$ is the $n$-th jumping time of the Poisson process $N_t^{\alpha}$, $T_{n} := \inf \{t : N_t^{\alpha} = n\}$}, where $(\alpha_t)_{t\ge 0}$ is a continuous and increasing function that describes the rate at which we use the sequence of observations $(Y_n)_{n\geq 0}$. We denote $(\alpha_t)_{t\ge 0}$ the \textit{intensity} of the process. On the one hand, and from a theoretical point of view, using more observations improve the algorithm  accuracy and convergence rate, so it may seem natural to use large values for $\alpha_t$. On the other hand, and in practice, observations can be costly and limited, so one would like to limit their use as much as possible. 

\begin{figure}\label{fig:Qgraph}
\centerline{\includegraphics[scale=0.4]{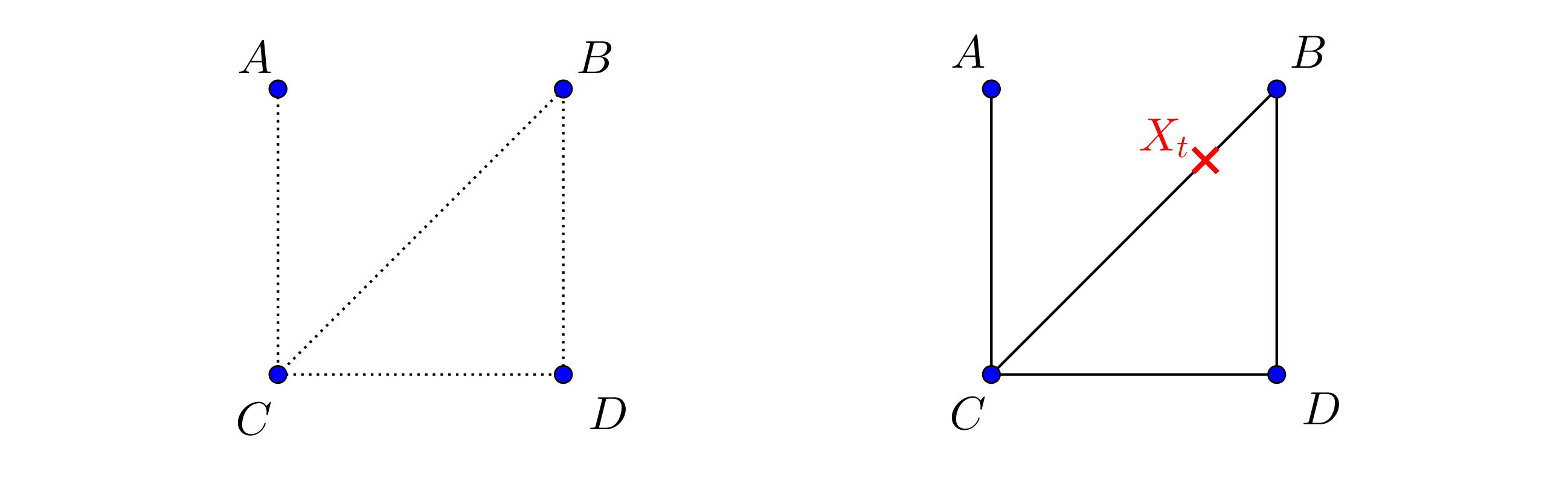}} 
\caption{\textbf{(left)} Example of discrete graph $G$, and \textbf{(right)} corresponding continuous version $\Gamma_G$. $X_t$ represents the current position of the algorithm in $\Gamma_G$. In the example its closest node in  $G$ is the node $B$.}
\end{figure}

We also emphasize that our algorithm runs on $\Gamma_G$, a continuous version of the discrete graph $G$, where each edge $e=(u,v)$ of length $L_e$ is seen as an interval $[0,L_e]$ such that an extremity of this segment corresponds to one of the nodes of the edge (see illustration Fig.~\ref{fig:Qgraph}). 
The process $X_t$ that represents the barycenter estimation at increasing times $t$, therefore lives on the graph edges and not just its nodes. Nevertheless, a current estimation of a central node can naturally be defined as the closest node to the position of $X_t$.

Note  that Alg.~\ref{algo:HLSA} is described for observations in $\{0, \cdots, K \}$. In an on-line context, the algorithm can then be used in three different ways: 
\begin{enumerate}
\item
If more than $K+1$ oservations are known: those used in the algorithm can be randomly picked-up.
\item
If less than $K+1$ observations are known and we won't have access to additional observations: Iterativally perform [(a) randomly shuffle the observations, and (b) use the shuffled obervations] until Alg.~\ref{algo:HLSA} ends. This strategy will be used in our tests.
\item
If less than $K+1$ observations are known and we will have access to additional observations: Use Alg.~\ref{algo:HLSA} on currently known observations to have a first guess of the graph barycenter. Wait for new observations to make the barycenter estimation more accurate with additional iterations of Alg.~\ref{algo:HLSA}.
\end{enumerate}

The key issue with this strategy on large graphs, which motivates the present paper, is that the deterministic move (row 7 of Alg.~\ref{algo:HLSA}) requires to compute the shortest path from $X_k$ to $Y_{k+1}$. To achieve this, we use of a stanard Dijkstra's algorithm which is particularly demanding in terms of computational times, especially when computed $K+1$ times. 
The solution of \cite{us} is then to pre-compute once for all the shortest distances between all node pairs and then to use this information for a quick algorithm execution. Computing these distances is $\#N$ times slower than computing the shortest path between two nodes, where $\#N$ is the number of nodes in $G=(N,E)$.
This solution then makes sense when  $K+1$ is larger than $\#N$, or when multiple runs of the algorithm will be performed on the same graph, \textit{e.g.} in order to evaluate the barycenters related to different observation sets $Y$. The major drawback of this strategy is however that it requires to store a $\#N \times \#N$ matrix in memory, which is unrealistic when $\#N$ is large. Moreover, the algorithmic cost of a Dijkstra's algorithm on our weighted graphs is anyway $\mathcal{O}(\#N^2)$ and therefore does not scale at all to large graphs. We then propose to use the multiscale solution described in the following section.

\section{Multiscale barycenter estimation}\label{sec:MultiscaleBE}

\subsection{General framework}\label{ssec:GenFramework}

Our method is motivated by a property of geometrical decomposition of the barycenter in Euclidean spaces. We will describe this property section~\ref{ssec:GraphPartition} and how we heuristically extend it to graph structures. In practice, an analysed graph $G$ will be parcelized (equivalently clusterized) into sub-graphs  and its barycenter will then be estimated using the \textit{Divide and Conquer} strategy given Alg.~\ref{alg:MBE}. This algorithm is directly related to Alg.~\ref{algo:HLSA} and the graph partition properties introduced section~\ref{ssec:GraphPartition}. Computation of the sub-graphs $G_i$ from a partition on large graphs is also discussed section~\ref{ssec:CptSubGraphs}.
Importantly, items 2 and 4 are however not as obvious as they may appear for two main reasons: (1) They should be scalable on large datasets, and (2) pertinent heuristics have to be used to define $\usg$ and $\mltg$, the simplified versions of $G$, so that they lead to accurate barycenter estimates of $G$. These two items are then presented sections \ref{ssec:SubGraphCpt} and \ref{ssec:MltGraph}.

\begin{algorithm}
\caption{Multiscale barycenter estimation}
\label{alg:MBE}
\begin{algorithmic}[1]
  \REQUIRE Graph $G=(N,E)$.
  \STATE Partition $G=(N,E)$ is partitioned into $I$ sub-graphs $G_i=(C_i,E_i)$.
  \STATE Undersample $G=(N,E)$ in $\usg=(\usv,\use)$, where each node of  $\usg$ represents a compact description of sub-graph $G_i$.
  \STATE Estimate the barycenter $\tilde{b}$ of $\usg$  using Alg.~\ref{algo:HLSA}.
  \STATE Compute a multiscale graph $\mltg=(\mltc,\mlte)$ with the nodes of $G$ in the subgraph of $\tilde{b}$ and the nodes of $\usg$ elsewhere.
  \STATE Estimate the barycenter $\bar{b}$ of $\mltg$  using Alg.~\ref{algo:HLSA}.
  \RETURN Node $\bar{b}$ estimated as the barycenter of $G$
\end{algorithmic}
\end{algorithm}

\subsection{Graph partition}\label{ssec:GraphPartition}

It is well known that for $n$ points, $(A_i)_{i=1,\cdots, n}$ of an affine space and an associated sequence of scalars $(a_i)_{i=1,\cdots,n}$ of non-null sum, the barycenter is defined at a single point $G=\mathrm{bar}\left((A_i,a_i)\right)_{i=1,\cdots,n}$ such that: 
\begin{equation*}
\sum_{i=1}^n a_i\bf{GA_i}=\bf{0}.
\end{equation*}
Suppose now that the nodes are partitioned into two node sets $I$ and $J$ and that the corresponding $(a_i)$ have a non-null sums. We then denote $G_I=\mathrm{bar}\left((A_i,a_i)\right)_{i\in I}$ and $G_J=\mathrm{bar}\left((A_i,a_i)\right)_{i\in J}$. The  decomposition property  states that the barycenter of the $n$ points is the barycenter of the two sub-barycenters, meaning:
\begin{equation}\label{prop:decompositionAffineBary}
\mathrm{bar}\left((A_i,a_i)\right)_{i=1\ldots n}=\mathrm{bar}\left( (G_I,\sum_{i\in I} a_i),(G_J,\sum_{i\in J}a_i)\right).
\end{equation}
This property can be iterated multiple times and still holds for $k$ partitions of this type, $k\le n$. 

Our multiscale graph barycenter estimation strategy is directly inspired by this property. 
In order to use a similar method on graphs, we use partitions (or clusters) $(C_i)_{i=1\ldots k}$ of the nodes set $N$ in $G$. We also define a sub-graph as:
\begin{definition}\label{def:associatedSubgraph}
For $C_i\subset N$, a subset of the nodes set, we call associated sub-graph $G_i$ a graph $G_i=(C_i,E_i)$ formed by all edges of the initial graph $G$, connecting two points of $C_i$. In other words, the edges set of $G_i$ is: 
\begin{equation}\label{def:E_i}
E_i=\{e=(e_-,e_+) \in E\ | e_-,e_+\in C_i\}.
\end{equation}
\end{definition}
A partition $\mathsf{P}=(C_i)_{i=1\ldots k}$ is called valid and can be used to compute graph barycenters, if it satisfies the following conditions:
\begin{enumerate}
\item The subsets $C_i$ are disjoint and their union contains all the nodes, \textit{i.e.} $N=\scup\limits_{i=1}^k C_i$;
\item The weight associated to each subset is non-null: $\forall 0\le i\le k$,  $\nu(C_i)\neq 0$;
\item Each part $C_i$ the associated sub-graph $G_i$ is connected.
\end{enumerate}
Note that condition (2) is implicit in our framework since $\nu$ charges each node, $\nu(y)>0, \ \forall y\in N$. Condition (3) is very important since the notion of barycenter in Alg.~\ref{alg:MBE} is only defined for connected graphs. 
Interestingly, a wide variety of established clustering algorithms efficiently define valid partitions of the nodes, even on large graphs \textit{e.g.} \cite{pythonAgglomClust}, that is based on \cite{NewmanPRE2004}.
We will therefore not develop this discussion in our paper and focus instead on the definition of   $\usg$ and $\mltg$.
First we define neighboring clusters w.r.t. $G$.
\begin{definition}\label{def:ngbClusters}
 Two disjoint subsets $C_i, C_j \subset N$  of the graph $G=(N,E)$ are neighboring clusters, denoted $C_i\sim C_j$, if there exists a pair of nodes $v_i\in C_i $ and $v_j\in C_j$ that are neighbors in $G$:
\begin{equation}
C_i\sim C_j \iff \exists\  v_i\in C_i, v_j\in C_j \mbox{ such that } (v_i,v_j)\in E.
\end{equation}
\end{definition}

In what follows, the information contained in the sub-graphs $G_i=(C_i,E_i)$ described above will be summarized in the graphs $\usg$ and $\mltg$ (the upscale and multi-scale versions of $G$). We can remark that the union of the edges in all clusters $C_i$ does not contain all the edges of the initial graph. The remaining edges will then be used to define the edges of  $\usg$ and $\mltg$.


\subsection{Computing the sub-graphs $G_i=(C_i,E_i)$}\label{ssec:CptSubGraphs}

The sub-graphs $G_i=(C_i,E_i)$ will be the key to subsample $G$ in $\usg$ and $\mltg$.
Here, we consider as known the partition $(C_i)_{i=1\cdots k}$ of the node set $N$. 
We then use Alg.~\ref{alg:SGC} to realistically compute the associated sub-graphs $G_i=(C_i,E_i)$.

\begin{algorithm}
\caption{Sub-graphs $G_i(C_i,E_i), \, i \in \{1,\cdots, k\}$ computation.}
\label{alg:SGC}
\begin{algorithmic}[1]
  \REQUIRE Graph $G=(N,E)$.
  \REQUIRE Nodes partition $(C_i)_{i=1\cdots k}$.
  \FOR{$i=1:k$}
    \STATE Create a void edge list $E_i$
    \STATE Create a void list $B_i$ of nodes information at the boundary of $G_i$ 
  \ENDFOR
  \FORALL{$e = (e_{-},e_{+}) \in E$}
    \IF{$e_{-}$ and $e_{+}$ are in the same cluster $C_i$}
      \STATE Add $e$ to $E_i$
    \ELSE
      \STATE We consider $e_{-}$ in cluster $i$ and $e_{+}$ in cluster $j$. 
      \STATE Add $[e_{-},e_{+},weight(e_{-},e_{+}),j]$ to $B_i$.
      \STATE Add $[e_{+},e_{-},weight(e_{+},e_{-}),i]$ to $B_j$.
    \ENDIF
  \ENDFOR
  \RETURN The sub-graphs $G_i(C_i,E_i), \, i \in \{1,\cdots, k\}$.
  \RETURN Sub-graphs boundary information $B_i, \, i \in \{1,\cdots, k\}$.
\end{algorithmic}
\end{algorithm}

At a first sight, the algorithmic cost of Alg.~\ref{alg:SGC} appears to be $\mathcal{O}(\#E)$. Checking the cluster of $e_{-}$ and $e_{+}$ (row 5)  however has an algorithmic cost $\mathcal{O}(\#N)$ if improperly coded. In our program, the node clusters are coded in Python dictionaries, making this task $\mathcal{O}(1)$ in average \cite{pythonTimeComp}. Using lower level programming languages, such as C++, the node identifiers could be first replaced by integers between $0$ and $\#N-1$ and then their clusters would be stored in a vector of size $\#N$, the cluster of node $i$ being stored at the $i^{th}$ entry of this vector. To be efficient this pre-treatment requires to sort the node labels; which has typically a cost $\mathcal{O}(\#N\log{(\#N)})$. This is for instance the case by using the standard C++ function std::sort. 
Once the node labels sorted, checking a node label is finally $\mathcal{O}(1)$, so this strategy also scales well to large graphs.

Finally, note that our algorithm not only computes the sub-graphs $G_i$ but also a compact information of their boundaries $B_i$. This boundary information will help us define a scalable strategy to generate the edges of $\usg$ and $\mltg$.

\subsection{Computing the subsampled graph $\usg$}\label{ssec:SubGraphCpt}

The definition of $\usg=(\usv,\use)$, the subsampled version of $G=(N,E)$, is only performed by using the sub-graphs $G_i=(C_i,E_i)$ and the boundary information $B_i$, $i \in \{1,\cdots, k\}$.  
Each node $\tilde{v}_i$ of $\usv$ indeed depends on the properties of $G_i$ and each edge of $\use$ depends on the $B_i$.

Every cluster $C_i$ is represented in $\usv$ by a single node $v_i$. The edge set of the new graph is defined by:
\begin{equation}
\use=\{(v_i,v_j) \mbox{ with } v_i,v_j\in\usv, v_i\in C_i, v_j\in C_i \mbox{ and } C_i\sim C_j\}.
\end{equation}
There exists an edge between two nodes if and only if their respective clusters are neighboring clusters and the length of each new edge is defined as the distance between its extremities in the subgraph $G_{ij}=(C_i\cup C_j, E_{ij})$, associated to $C_i\cup C_j$.

\begin{figure}[htb!]
\begin{center}
\includegraphics[width=0.90\linewidth]{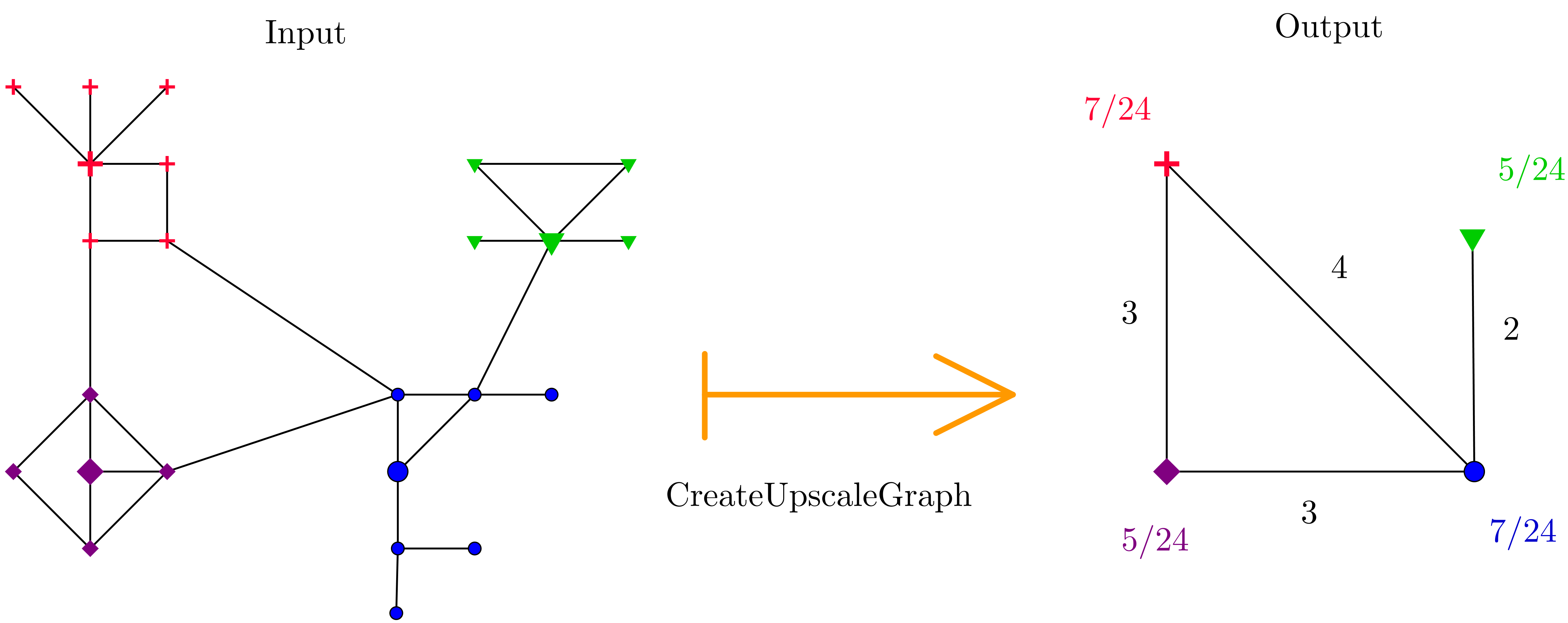}
\caption{
On the left hand side we have an initial graph $G$ partitioned in $4$ communities (the nodes of each community are represented by different symbols and colors). We consider $\nu$ as the uniform probability on its nodes set and all edges have length $1$. We take a subset $\usv$ formed of nodes represented in a larger size. The upscale graph $\usg$, is represented on the right hand side, along with the length of the new edges and the new probability corresponding to each node. We have chosen to represent the distribution $\nu_{\usg}$, instead of a sequence $(\usy_n)_{n\ge 0}$, because it is easier to visualize.  }
\label{fig:CUG}
\end{center}
\end{figure}

From a mathematical point of view, the probability associated to each node is the total probability of the cluster that contains it:
\begin{equation}\label{def:probaUpscale}
\forall v\in\usv, \quad \nu_{\usg}(v)=\nu(C_i),\quad \mbox{ where } C_i \mbox{ is such that } v\in C_i.
\end{equation}
The definition of the associated probability measure for the upscale graph in (\ref{def:probaUpscale}) is the analog of summing the scalars in the affine case in (\ref{prop:decompositionAffineBary}). 
Accessing independent random variables $(Y_n)_{n\ge 0}$ distributed according to $\nu$, we can easily define another sequence $(\usy_n)_{n\ge 0}$ of i.i.d. random variable of law $\nu_{\usg}$:
\begin{equation}\label{def:upscaleY}
\usy_n=c_i\in \usv\cap C_i \quad \mbox{ if and only if }\quad Y_n\in C_i.
\end{equation} 
From the simulation point of view, when we have access to $(Y_n)_{n\ge 0}$, \eqref{def:upscaleY} means that every time a node in a cluster $C_i$ is given by the sequence, we see it as the unique node $c_i$ that represents the cluster in the upscale graph. 
An example of this procedure is illustrated on a simple graph in Figure \ref{fig:CUG} and further details are presented in the following subsections.  \\

\subsubsection{Computing the subsampled nodes $\usv$}\label{sss:subnodes}

A natural strategy to compute each subsampled node $\tilde{v}_i$ is to define it as the barycenter of $G_i(C_i,E_i)$ using Alg.~\ref{algo:HLSA}. This however requires to have a sufficient number of observations $Y$ in $C_i$. We then instead randomly draw a node $C_i$ to define $\tilde{v}_i$, which has also the advantage of having a negligible computational cost. 

\subsubsection{Computing the distance between $\tilde{v}_i$ and the  boundaries of $G_i$}\label{sss:computing borders}

Defining the distances between the nodes $\usv$ in $\usg$ is the trickiest step of our strategy. Consider two neighbor sub-graphs $G_i$ and $G_j$. A quick and simple strategy would be to define the distance between $\tilde{v}_i$ and $\tilde{v}_j$ as equal to the diameter of $G_i$ plus the diameter of $G_j$ divided by two. Early tests performed using this strategy have however not been satisfactory. We will then define in next subsection this distance as the shortest distance between $\tilde{v}_i$ and $\tilde{v}_j$ in the union of the sub-graphs of $G_i$ and $G_j$, plus the edges of $G$ linking these subgraphs that are saved in $B_i$ and $B_j$. 

For each region $i$, we compute the distance between $\tilde{v}_i$ and all boundary nodes of sub-graph $G_i$, \textit{i.e.} the $B_i[j][1]$, $j \in \{1,\cdots,\#B_i\}$ (see rows 9 to 11 of Alg.~\ref{alg:SGC}). This can be done by running a Dijkstra's algorithm as discussed in the end of Section~ \label{sec:SA} for the whole graph, where it was too costly. A fundamental remark here is that although this strategy was far too computationally consuming for the whole graph (algorithmic cost is $\mathcal{O}(\#N^2)$), it becomes realistic on much smaller clusters $G_i(C_i,E_i)$. It can also be straightforwardly parallelized on different sub-graphs. After having computed the distances, we add them to the corresponding sub-lists of $B_i$. Each sub-list $j$ of the list $B_i$ has then the following structure:
\begin{itemize}
\item
$B_i[j][1]$: Node of $G_i$ at the cluster boundary.
\item
$B_i[j][2]$: Node outside of $G_i$ linked to $B_i[j][1]$.
\item
$B_i[j][3]$: Distance between $B_i[j][1]$ and $B_i[j][2]$.
\item
$B_i[j][4]$: Cluster of $B_i[j][2]$.
\item
$B_i[j][5]$: Distance between $B_i[j][1]$ and $\tilde{v}_i$.
\end{itemize}

\subsubsection{Computing the subsampled edges $\use$}

The definition of the edges $\use$ only depends on the boundary information $B_i, \, i \in \{1,\cdots, k\}$, as explained Alg.~\ref{alg:cptSBE}. 
\begin{algorithm}
\caption{Subsampled graph edges $\use$ computation.}
\label{alg:cptSBE}
\begin{algorithmic}[1]
  \REQUIRE Clusters boundary information $B_i, \, i \in \{1,\cdots, k\}$.
  \FOR{$i=1:k$}
    \FOR{$j=1:\#B_i$}
      \IF{$i<B_i[j][4]$}  \label{line:undirected}
        \STATE $\bar{i}=B_i[j][4]$
        \STATE Identify $\bar{j}$ so that $B_i[j][1]==B_{\bar{i}}[\bar{j}][2]$ and $B_i[j][2]==B_{\bar{i}}[\bar{j}][1]$.
        \STATE Compute $TmpDist=B_i[j][3]+B_i[j][5]+B_{\bar{i}}[\bar{j}][5]$.
        \IF{$\use$ does not contain the edge ($\tilde{v}_i$,$\tilde{v}_{\bar{i}}$) or its distance is $>TmpDist$}
          \STATE Add or update edge ($\tilde{v}_i$,$\tilde{v}_{\bar{i}}$) to $\use$ with distance $TmpDist$.
        \ENDIF
      \ENDIF
    \ENDFOR
  \ENDFOR
  \RETURN Subsampled graph edges $\use$.
\end{algorithmic}
\end{algorithm}

Remark that the test line~\ref{line:undirected}, of Alg.~\ref{alg:cptSBE} is performed to avoid having multiple edges linking the same nodes as the graphs are undirected.
This algorithm is again computationally reasonable as the main double for loop first depends the number of clusters  and the number of nodes at the clusters boundaries. The instructions in this double loop are also reasonable as they are linearly related to limited number of edges and nodes.

\subsubsection{Projecting the observations $Y$ from the nodes of $G$ to those of $\usg$}

We recall that efficient techniques were described in section~\ref{ssec:CptSubGraphs} to find the sub-graph $G_i$ associated to each node of $N$. We use the same technique to project the obervations $Y$ on $N$ to each node $\tilde{v}_i$ of $\usg$. The node $\tilde{v}_i$ indeed represents all the nodes $C_i$ of $G_i$.

\subsection{Multiscale graph}\label{ssec:MltGraph}

\subsubsection{Motivation}


A straightforward extension of the decomposition of the barycenter in the Euclidean case, to the context of graphs, can be described as follows: take a valid partition, compute the barycenter of each cluster, create a new subgraph (as explained in subsection \ref{ssec:SubGraphCpt}) and finally compute its barycenter. This procedure induces thus a notion of centrality on the set of clusters.  

Therefore, by choosing each subsampled node $\tilde{v}_i$ as a Fr\'echet mean of $C_i$ (in subsection \ref{sss:subnodes}), creating the corresponding subsampled graph $\usg$ as described above, and then estimating its barycenter, in a sense, we obtain a central cluster.  
If the chosen partition has a specific meaning, this procedure can have an interest on its own, allowing us to study some larger scale properties of the graph. For example if each cluster $C_i$ represents a community, this is a natural way of defining a central community.

Independently of the method used to define the starting points $\tilde{N}$ (randomly chosen representatives of each subset $C_i$ or  estimated barycenter), since the graph does not have the same properties as the euclidean space, a barycenter $\tilde{b}$ of the subsampled graph $\usg$ is not necessarily a barycenter of the initial graph. However, for reasonable partitions, one might expected the Fr\'echet mean of the initial graph not to be far from the central community that contains $\tilde{b}$. This assumption motivates the next step in our approach: building the multiscale graph as detailed hereafter.

\subsubsection{Definition of a multiscale graph}
For a valid partition of the nodes set $\mathsf{P}=(C_i)_{i \le k}$, let $(G_i,\nu_i)_{i\le k}$ denote the associate sub-graphs with their respective probabilities measures, defined in subsection \ref{ssec:CptSubGraphs}. With the notations introduced in \ref{ssec:SubGraphCpt}, let $\usg=(\usv,\use)$ be an up-scale version of $G$ corresponding to the partition $\mathsf{P}$. In what follows we define $\mltg=(\mltv,\mlte)$ the \textit{multi-scale version } of $G$ with respect to $(\usg,C)$, where $C$ is an element of $\mathsf{P}$, and $\mltn$ the corresponding probability measure. The definition of the nodes set and the associated probability are straightforward. $\mltv$ contains the nodes of $\usg$ and $C$:
\begin{equation}\label{def:multivertices}
\mltv= \usv\cup C,
\end{equation}
and $\mltn$ redistributes the mass of $C$ to its nodes, while leaving the others values of $\nu_{\usg}$ unchanged:
\begin{equation}\label{def:mlt:proba}
\mltn(v)=\begin{cases} \nu(v) & \mbox{if } v\in C\\
\nu_{\usg}(v) &\mbox{if } v \in \usv
\end{cases}
\end{equation}
The edge set $\mlte$ contains the edges of $\usg$, except those that were added to $c$ , the node that represents the cluster $C$ in the up-scale graph, and all internal edges of $C$. On top of that we add new edges going from boundary of $C$ to the nodes corresponding to its neighboring clusters:

 \begin{equation}
 \mathrm{BorderEdges}(C,\mathrm{P})=\{(v,c_i)|\ v\in C, \exists v_i \in C_i \mbox{ with }(v,c_i)\in E\}
\end{equation}  
The length of such an edge is defined as the initial distance between its extremities in the subgraph corresponding to $C\cup C_i$. Now the set of edges $\mlte$ can be written as:
\begin{equation}\label{def:multi:edges}
\mlte=E_C\cup \left(\use\setminus \{e| e\sim c\}\right)\cup \mathrm{BorderEdges}(C,\mathrm{P}),
\end{equation} 
where $e\sim c$ means that $c$ is a node of $e$.

\subsubsection{Computing the multiscale graph $\mltg$}\label{ssec:MltGraphCpt}
As explained Alg.~\ref{alg:MBE}, the multiscale graph $\mltg=(\mltc,\mlte)$ has the nodes of $G$ in the subgraph of $\tilde{b}$ and the nodes of $\usg$ elsewhere, where $\tilde{b}$ is the estimated the barycenter of $\usg$  using Alg.~\ref{algo:HLSA}.

We denote $\tilde{i}$ the label of the sub-graph containing $\tilde{b}$. The construction of $\mltg$ is then the same as the one of $\usg$, except in the sub-graph $G_{\tilde{i}}$, where the original nodes and edges are preserved. At the boundary between $G_{\tilde{i}}$ and the subsampled domain the distance given to the edge is slightly different to row 6 of Alg.~\ref{alg:cptSBE}. The distance between a sub-graph representative and the sub-graph boundary (\textit{i.e.} $B_i[j][5]$ or $B_{\bar{i}}[\bar{j}][5]$) in $G_{\tilde{i}}$ is obviously not considered.\footnote{The nodes on the central graph are not connected to its boundary ? }


A fine representation of $G$ is then constructed in  $G_{\tilde{i}}$, the central cluster, and a coarse representation of $G$ is constructed elsewhere. The goal of this multiscale graph is to make it possible to finely  estimate the barycenter of $G$ with reasonable computational resources.

Of course, from a simulation point of view, constructing a sequence of random variables distributed according to $\mltn$, the probability distribution defined in \eqref{def:mlt:proba} is straightforward once we have access to $(Y_n)_{n\ge 0}$ of law $\nu$. We simply set:

\begin{equation}\label{def:multi:y}
\hat{Y}_n=\begin{cases}Y_n, &\mbox{ if }Y_n\in C\\
v_i, &\mbox{ if }Y_n\in C_i.
\end{cases}
\end{equation}

\subsection{Illustration on the Parisian subway}\label{ssec:Illustration}

\begin{figure}[htb!]
\begin{center}
\includegraphics[width=0.70\linewidth]{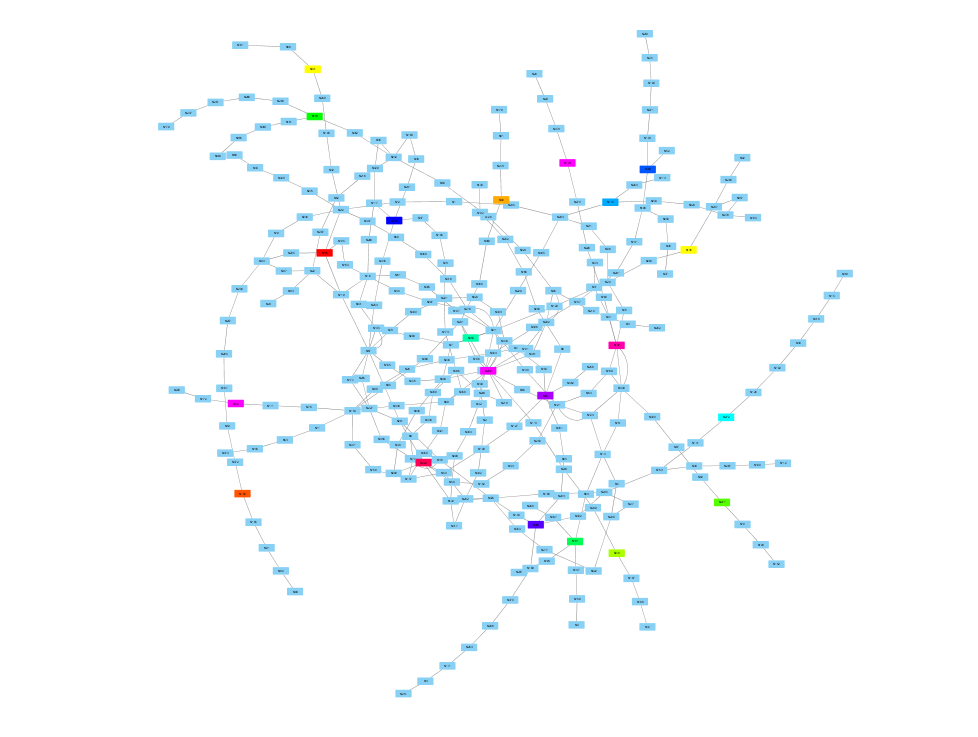}
\caption{
Complete graph of Parisian Metro $G$. The colored nodes were randomly drawn in the precomputed sub-graphs $G_i$ and are be the nodes $\usv$ of $\usg$ (see Fig.~\ref{fig:ReScaledMetro}).
}
\label{fig:completeMetro}
\end{center}
\end{figure}

We now illustrate what can be $G$, $\usg$ and $\mltg$ on the Parisian subway network, which we already used in \cite{us}. Fig.~\ref{fig:completeMetro} represents the complete Parisian subway network. The graph was downloaded at \url{http://perso.esiee.fr/~coustyj/EnglishMorphoGraph/PS3.html}, has 296 nodes and 353 nodes. The nodes obviously represent the metro stations. Each edge is a connection between two stations and its  length is the time needed to go from one station to the other. 

\begin{figure}[htb!]
\begin{center}
\includegraphics[width=0.70\linewidth]{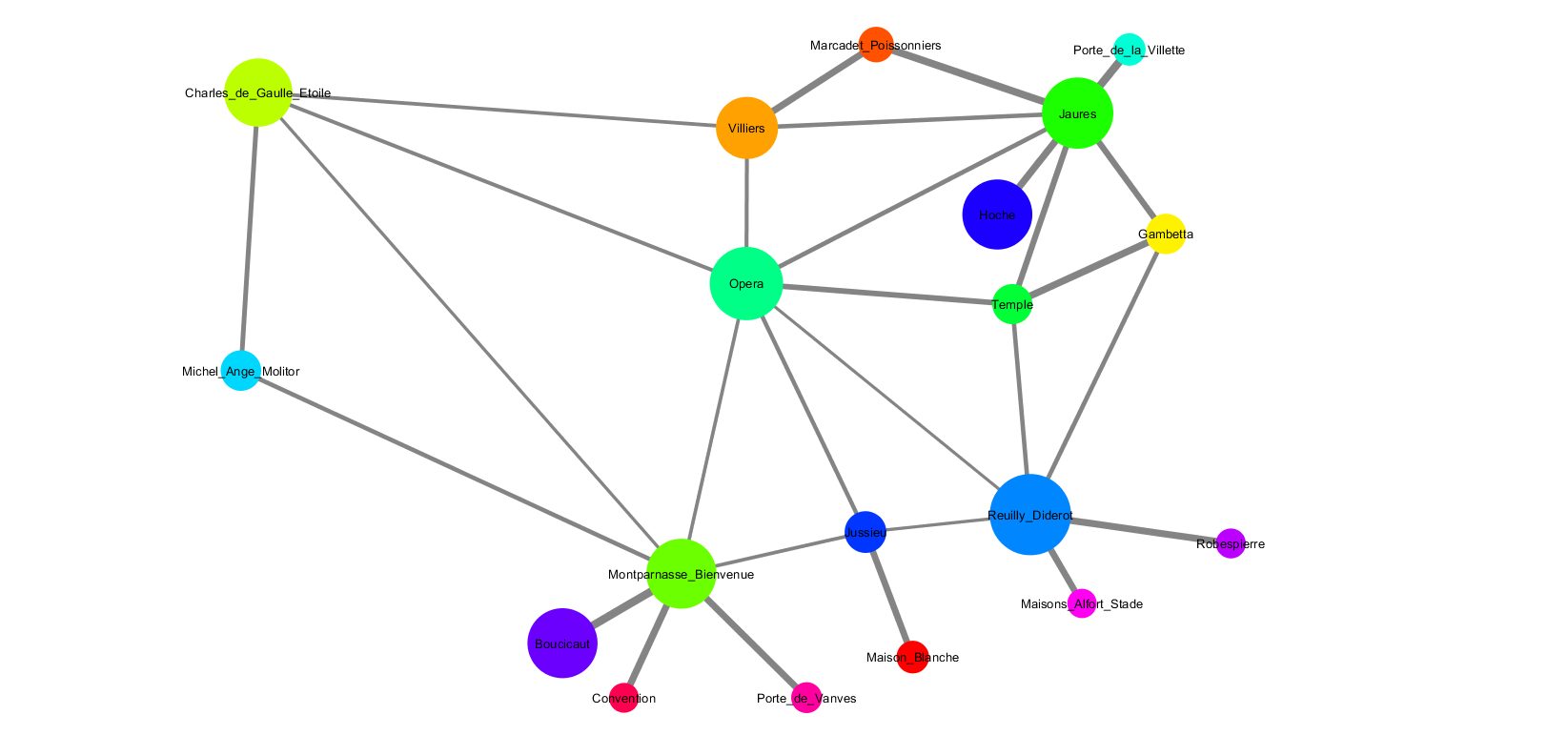}
\includegraphics[width=0.70\linewidth]{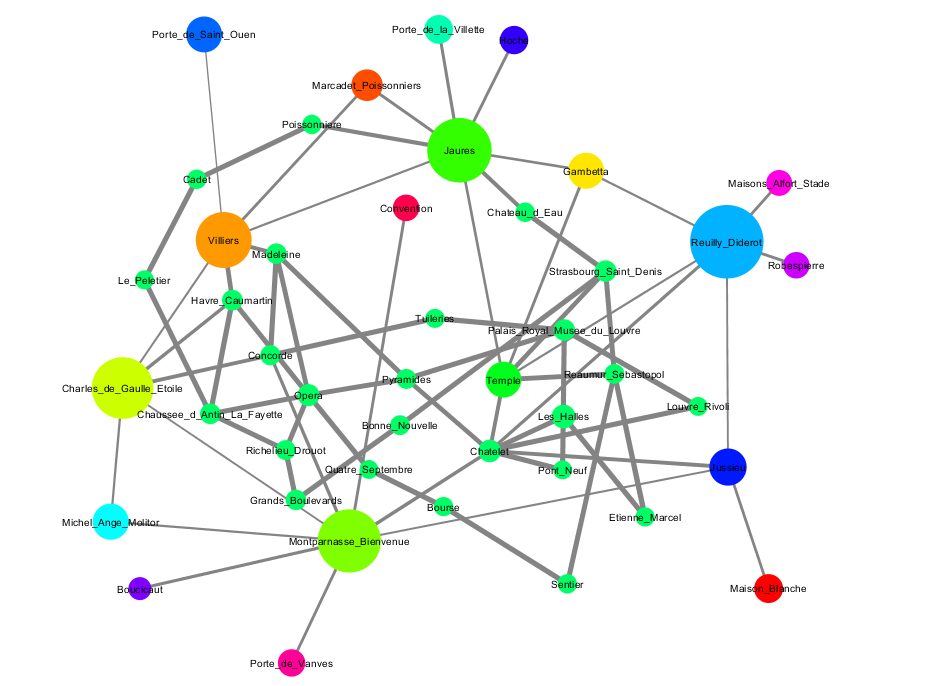}

\caption{
\textbf{(top)} Subsampled Parisian Metro graph $\usg$, and \textbf{(bottom)} multiscale Parisian Metro graph $\mltg$. The width of the edges is inversely proportional to the time needed to go from one station to the other.
}
\label{fig:ReScaledMetro}
\end{center}
\end{figure}

In Fig.~\ref{fig:ReScaledMetro}(top) we represent the subsampled graph $\usg$. One can see that the actual barycenter of the initial graph $G$ (Chatelet) is not included in $\usg$, and thus can't be estimated as its center. It can however be estimated in the multiscale graph $\mltg$ that is show Fig.~\ref{fig:ReScaledMetro}(bottom).

\medskip

\section{Results}

\subsection{Results on small graphs}
In order to validate our strategy, we tested it on three small graphs for which we have access to the ground-truth barycenter. The first one is the Parisian metro network descriebed in \ref{ssec:Illustration}. The other two subgraphs of Facebook from the Stanford Large Network Dataset Collection\footnote{\url{https://snap.stanford.edu/data/}}: \textbf{(FB2000)} has 2000 nodes and 37645 edges and \textbf{(FB4000)} has 4039 nodes and 88234 edges and fully contains (FB2000).

  We performed three type of tests with default parameters: 
\begin{itemize}
\labitem{1}{ss} Single scale estimation using Alg.\ref{algo:HLSA}.
\labitem{2}{ms} Multi-scale estimation; in the upscale graph each cluster is represented by its barycenters (estimated using Alg.\ref{algo:HLSA}).
\labitem{3}{mmr} Multi-scale estimation; in the upscale graph each cluster is represented by a node sampled at random (uniformly) among its nodes. 
\end{itemize} 
We ran $100$ Monte Carlo simulations for each strategy. A run is considered successful if the returned node is the true barycenter of the graph. We sum up the results in Table \ref{res:FB2000}.

\begin{table}

\begin{center}
\begin{tabular}{|c|c|c|c|}
\hline
 & Single scale & Multi-scale & Multi-scale random\\
 \hline
 Paris Network & 100\% & 97 \% & 97 \% \\
 \hline
FB2000 & 100\% & 100\% & 100\%\\
\hline
FB4031 & 100\% & 80\% & 73 \% \\
\hline
\end{tabular}
\caption{Succes ratio obtained with each strategy on 100 Monte Carlo runs for the Parisian metro network and FB2000 graph. \label{res:FB2000}}
\end{center}

\end{table}
As one can see in Table \ref{res:FB2000}, on the first two graphs, Parisian metro network and the FB2000, the performance of the algorithm does not seem influenced by choice of the nodes that represent each cluster in  the upscale graph. On the third one, FB4000, the success ratio decreases slightly, but the algorithm still performs rather well when the representative nodes in the preliminary phase are chosen at random. 
This is what motivated us to apply this strategy on larger graphs in order to reduce computational cost. We do not claim however that this is efficient in any framework. Indeed, depending on the graph's structure and the initial partition, there exist cases where this first approximation is crucial and the choice of the representative node of a cluster can directly impact the quality of the results.

%
%

\subsection{New York Urban Area}
\begin{figure}[htb!]
\begin{center}
\includegraphics[width=0.70\linewidth]{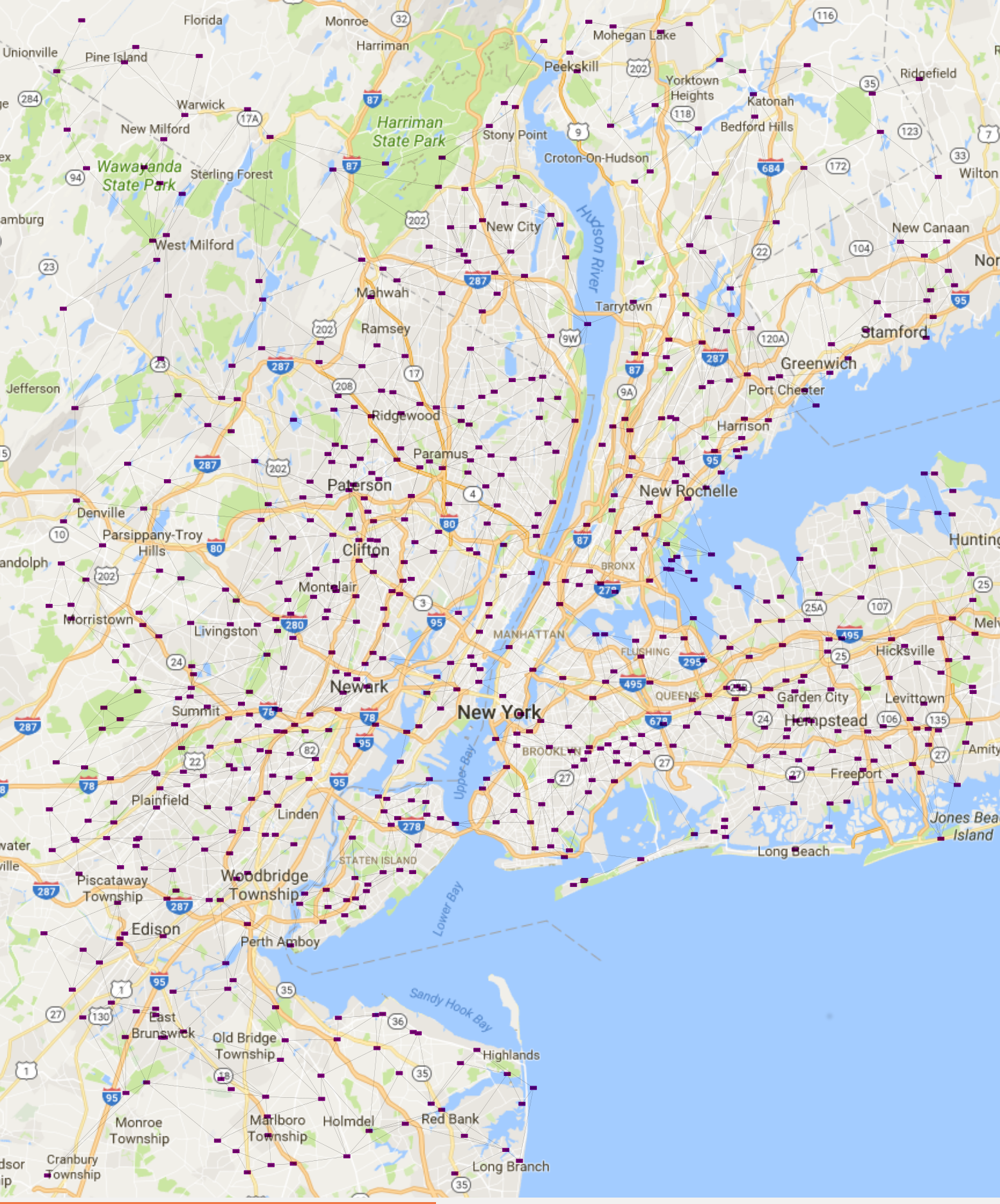} 
\caption{New York urban area. Image obtained using Cytoscape and \copyright Google Maps. Purple points represent estimated centers of the NG700 partition.}
\label{fig:NY:Complete}
\end{center}
\end{figure}

After measuring the stability and accuracy of our method on small graphs, where the ground-truth barycenter is known, we have chosen to test it on a graph formed by the the crossroads in a rather large New York urban area. We referred to it as such by convenience, but the area is not limited to the state of New York, see Figure \ref{fig:NY:Complete}. The graph has 264.346 nodes and 733.846 edges and can be found on the website of the Center for Discrete Mathematics and Theoretical Computer Science \footnote{\url{http://www.dis.uniroma1.it/challenge9/download.shtml}}. On the website it is mentioned that some gaps might exist and thus the graph does not necessarily contain all crossroads. The nodes are the GPS coordinates of crossroads and the edges represent the streets between them. The distance considered between two nodes is the physical one, and not the transit time. Furthermore, the graph is undirected, namely each street allows travel in both directions. \\

We have performed two types of preliminary clustering. One, based on a bottom up approach, meant to provide clusters of homogeneous size, and another, based on a Markov Clustering algorithm developed by Stijn van Dongen \cite{mcl}. The graph partition obtained with the first method has 700 clusters (from now on it will be referred to as GP700) and the second one has 1776 (we will referred to it as MCL12).\\

\paragraph{Technical details}In terms of memory cost, this kind of graph can definitely not be handled with the method proposed in \cite{us}. A rough estimation suggests that the associated matrix distance would need around 360 GB of memory, whereas this new method employs far less resources, being of the order of 15GB (or less). The computational time is a bit long, but reasonable. Using the default parameters it takes, in average, 3 h 30 min  for the GP700 and 7 h 30 min for the MCL12. This can be easily improved by paralelizing the computation of the subsampled graph $\usg$, namely the computation of its sub-sampled nodes (see sub-section \ref{sss:subnodes}) and the informations related to the borders (see subsection \ref{sss:computing borders}), needed for the computation of its edges (in Alg.\ref{alg:cptSBE}). It is not surprising that the barycenter's estimation on GP700 is faster, since the clusters have a more homogeneous size and are thus easier to handle.\\


Since the graph is too dense to visualize, we have chosen to use the GP700 partition in order to facilitate Figure \ref{fig:NY:Complete}. To be more precise, we have used the upper scale approximation procedure described in Section \ref{ssec:SubGraphCpt} to form a new graph from the estimated centers of each cluster. A visualization using the GPS coordinates of the nodes was created with the aid of the Cytoscape software. This illustration was afterwards overlaid on a map of the area provided by \copyright Google Maps. The result is shown in Figure \ref{fig:NY:Complete}. The purpose of this figure is to give an idea of the area covered by the complete graph and not to show the exact position of each node in the upper scale graph. 

\subsubsection{Results on two different partitions}
\begin{figure}
\begin{center}
\includegraphics[width=10cm]{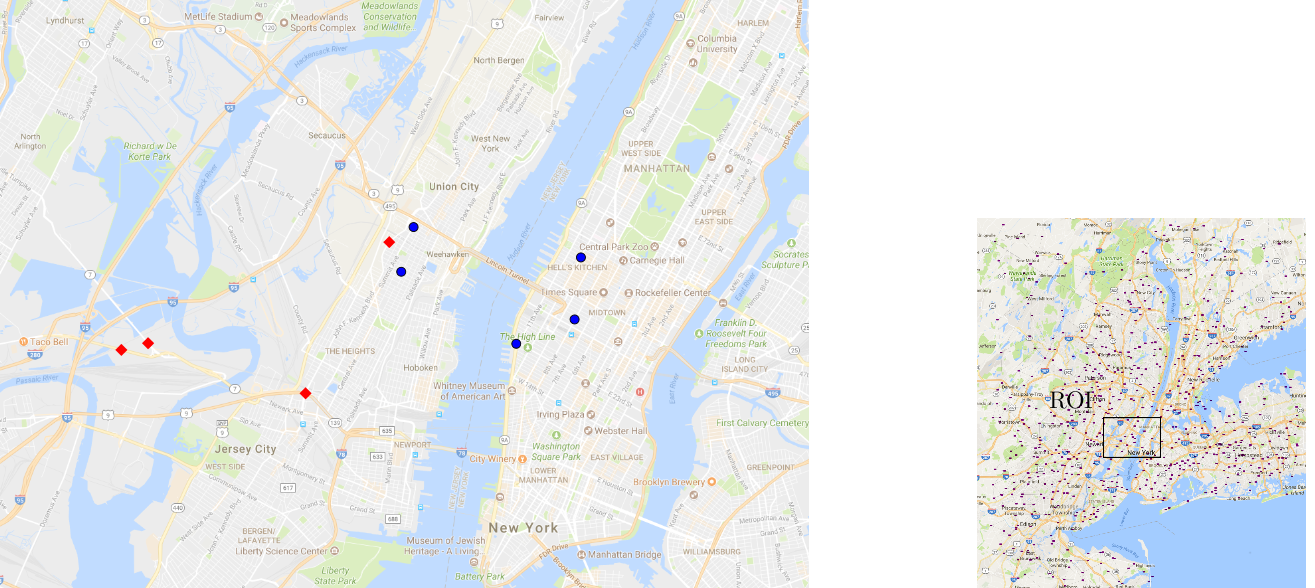}
\hfill
\caption{On the left hand side we have a general view over the complete New York graph represented in Figure \ref{fig:NY:Complete}. On the right hand side, we have the Region Of Interest (ROI). The  black round points are the barycenters obtained on the GP700 partition and the red diamond-shaped points are the barycenters obtained on MCL12. The right hand-side image was obtained using \copyright Geogebra and \copyright Google Maps.}
\label{fig:MCL12+GP700}
\end{center}
\end{figure}

We illustrate the results obtained in $4$ Monte Carlo runs on the partitions GP700 and MCL12  in Figure \ref{fig:MCL12+GP700}. Since we have the GPS coordinates of each node, we used them to represent the estimated barycenters with \copyright Google Maps. We computed the mean distance $MD$ between each set of barycenters directly on the graph:
$$MD(GP700)= 35 \quad \quad MD(MCL12)=50. $$

%
%
\subsubsection{Parameters influence}
In this section we illustrate the influence of the parameters on the results obtained on GP700. As explained in Section \ref{sec:SA}, the main parameters are: the temperature schedule $\beta$, the rate at which the observations are used $\alpha$ and the stopping time $T$.  We consider as default parameters, and denote $\beta^*, \alpha^*$ and $T_{\max}^*$, the parameters introduced in \cite{us}. 
\paragraph*{Influence of the temperature schedule.}
The temperature schedule $\beta_t$ is a very important parameter linked to simulated annealing. Large values of $\beta$ increase the convergence speed of the algorithm. However, if its value is too large, if it crosses a certain threshold that depends on the graph's structure, the algorithm might converge to a local minimum instead of a global one. The convergence of the simulated annealing is guaranteed from a theoretical point of view for logarithmic evolutions of the temperature( $beta_t=\beta \log t$), and it is with this type of growth that we established the value of the constant $\beta^*$. However, in practice a linear growth is more commonly used ($\beta_t=\beta t$) we have tested both versions. We run our algorithm five times for each set of parameters and computed the mean distance between the estimated barycenters in order to measure the stability.  As one can see on Table \ref{res:beta}, for a logarithmic schedule, increasing the constant $\beta$, reduces the variations of the results (because the algorithm converges faster). However for a linear evolution, increasing the value of the constant $\beta$ destabilizes the algorithm (probably because the algorithm tends to converge to local minimums).  

\begin{table}
\begin{center}
\begin{tabular}{|c|c|c|}
\hline
  Parameter &Estimated Center Log & Random center Linear  \\

  \hline
default & 55  & 35 \\
\hline
  $\beta_{0.25}$ & 92&  32  \\
  \hline
  $\beta_{0.5}$ & 88&  21 \\
\hline  
  $\beta_2$ & 44 &  45  \\
\hline  
  $\beta_4$ & 44 &  56  \\
 \hline
 \end{tabular}
\caption{Mean distance between the final centers obtained in $5$ runs for each set of parameters. \label{res:beta}}
\end{center}
\end{table}

\paragraph*{Influence of $\alpha$ and the stopping time. } In our tests, the rate at which we use the observations is calibrated with respect to the stopping time. The stopping time, for a barycenter estimation using Algorithm \ref{algo:HLSA}, is chosen as a function of the number of nodes in the graph: $$T^*_{\max}=0.1 \#V+100.$$  
The use of current observations is then distributed at a rate $\alpha_t$ that insures that between $T^*_{\max}$ and $T^{\star}_{\max}$ we use approximately $S^*=1000$ observations. In theory, a balance between the intensity $\alpha_t$ and the temperature $\beta_t$ is mandatory for the convergence of the simulated annealing. So increasing the stopping time, without increasing accordingly $S^*$, reduces the intensity rate $\alpha$.  As shown in Table \ref{res:alpha} this can be problematic, especially when we use a linear growth temperature.

\begin{table}

\begin{center}
\begin{tabular}{|c|c|c|}
\hline
  Parameter &Estimated Center Log & Random Center Linear  \\
\hline  
 
default & 55 &    35 \\
\hline

  $0.25T^*_{\max}$ & 68 &  33  \\
\hline 
  $0.5T^*_{\max}$ & 70 &  35 \\
  \hline
$2T^*_{\max}$ & 33 &  60  \\
  \hline
  $4T^*_{\max}$ & 54 & 104 \\
 \hline
 \end{tabular}
\caption{Mean distance (divided by $10^3$ to be easier to read) between the final centers obtained in $5$ runs for each set of parameters.\label{res:alpha} }
\end{center}
\end{table}
\paragraph*{Influence of the observations number.} As explained in Section \ref{sec:SA}, when we have access to a limited number of observations, there are multiple strategies available. For the current tests, we used strategy $2$, namely we supposed that less than the necessary number of observations are available, and thus once the list is exhausted, we shuffle it and reuse it. This is equivalent to accurately estimating the barycenter of the multiscale graph, endowed with the probability measure $\hat{\nu}$, corresponding to the known observations. 

For our tests, we generated a list of $O^*=10000$ observations by choosing nodes uniformly at random among the vertex set. We then created other lists with more or less observations to measure the influence of this number on our method. The results are summarised in Table \ref{res:obs}.  For the linear temperature schedule, the variance seems to be quite stable with respect to the number of observations. The fact that for a logarithmic evolution, the variance increases with the number of observations is not surprising since, less observations imply a more concentrated probability measure $\hat{\nu}$, and thus its barycenter might be easier to estimate. However, the stability of the algorithm should not be regarded as the ultimate guarantee of the quality of its results. The more observations we use, the more we get closer to estimating the barycenter that corresponds to the uniform probability measure on the entire graph. And even though this bias is not very important for a uniform measure we expected it to be more prominent on heterogeneous probability measures.

\begin{table}

\begin{center}
\begin{tabular}{|c|c|c|}
\hline
  Parameter &Estimated Center Log & Random Center Linear  \\
\hline  
 
default & 55 &   35 \\

  \hline
$ 0.01 O^*$ & 13 &  62 \\
\hline
 $0.1 O^* $&  58  &  81\\
\hline
 $10 O^*$& 88 & 58 \\
\hline
  $100 O^*$ &  82 & 44\\
\hline
\end{tabular}
\caption{Mean distance (divided by $10^3$ to be easier to read) between the final centers obtained in $5$ runs for each set of parameters.\label{res:obs} }
\end{center}
\end{table}



%

\subsection{Social network}
Finally we tested our method on a Youtube sub-graph downloaded from  the Standard Large Network Dataset Collection, of $1,134,890$ nodes and $2,987,624$ edges. Each node represents a user and the edges represent user-to-user links. Of course, the data is anonymized. Each edge is of length $1$ and the observations are uniformly sampled from the vertex set.

 \paragraph{Results} Using the default parameters, we ran the algorithm $4$ times and we obtained two different estimations of the barycenter: nodes $'1072'$ and $'663931'$.
The two nodes are not directly connected  and the distance between them is equal to $2$. We could say that this distance represents an average closeness, since it is slightly lower than the mean distance between each of them and all other nodes of the graph, which is approximately $2.95$ for   $ '1072'$ and $3.51$ for $'663931'$.  We do not know the ground-truth barycenter for this graph, but the results are quite stable and thus promising. The computation time for one run of the algorithm was around $64$ hours. The results are summarized in Table \ref{res:you}.
\begin{table}

\begin{center}
\begin{tabular}{|c|c|c|}
\hline
  Node & Frequency & Mean Distance\\
 \hline
 1072  & 3 & 2.95  \\
 \hline
663931 & 1 & 3.5 \\

\hline
\end{tabular}
\caption{The second column represents the number of times the node was select as a barycenter by our Algorithm on the $4$ Monte Carlo runs. The mean distance is the average distance from the selected barycenter to all other nodes of the network.}\label{res:you}
\end{center}

\end{table}
 
 \subsection{Conclusion}
\paragraph*{Memory cost.}$ $
 From a computational point of view, the multiscale approach drastically reduces the memory costs and thus can be used on larger graphs.

\paragraph*{Computational time.}$ $ Currently, the computational time is rather long, but as mentioned before it can be reduced by parallelizing some of the intermediate procedures. Moreover, the most costly part of the algorithm is the conception of the subsampled graph and this step does not need to be done at the arrival of a new observation in the online context.

\paragraph*{Online update of the barycenter.}$ $
As mentioned before, the most time costly operation is the creation of the upscale graph. However, the actualization of the barycenter on the multiscale graph at the arrival of a new observation is instantaneous. If the informations regarding the upscale graph are stored, we could even reset the algorithm this stage and update the  barycenter estimation on the subsampled graph. Assuming that the estimated central cluster would change with this new observation, creating a new multiscale graph and estimating a new barycenter on it (using the default parameters and the corresponding number of observations) would take less than $1$ minute for the New York graph (with the partitions we used for our test) and around $7$ hours for the Youtube graph. The time needed for this operation  depends a lot on the size of the clusters and not only on the size of the initial graph. For example, in our tests for the Youtube graph, computing the distances on the multiscale graph takes around $6$ hours and estimating its barycenter only one. 


\appendix{}

\section{Package description}

The \textit{LGC\_estim} package contains the strategy described in this paper. It is entirely written in Python and was tested using Python 2.7 and 3.8 \textcolor{red}{(Ioana, did you used Python 3.x? If yes, which version?: answer 3.4)}. Outside of Python modules that can be considered as standard (\textit{Numpy}, \textit{sys}, \textit{os}, $\cdots$), the only specific dependence of our package is \textit{NetworkX}. This module is widely used for graphs management and analysis in Python\footnote{For installation, please go to \url{https://networkx.github.io/documentation/networkx-1.11/install.html}}. Note that all our tests were made using the version 1.11 of NetworkX.

There are two ways to use the \textit{LGC\_estim} package:  It was primarly designed to be used as a script but it can be alternatively used as a  Python module. A \textit{README} file at the package root directory explains how to use it in both cases through simple examples. Data files representing the Parisian subway network of section~\ref{ssec:Illustration} are included in the \textit{data} directory to run the examples. Note that by simply executing the
command line \textit{python LargeGraphCenterEstimator.py}, as for any Python script, a help message will give intructions to follow to properly estimate graph barycenters. Note finally that using the \textit{LGC\_estim} package as a Python module requires to understand the key classes and functions we used in our code, but is doable as shown in the \textit{README} file.

\begin{acknowledgements}
The authors thank S\'ebastien Gadat and Laurent Miclo for insightful discussions.
\end{acknowledgements}


\bibliographystyle{spmpsci}      
\bibliography{PaperMPC}   


\end{document}